\documentclass{amsart}
\usepackage{amsthm}
\usepackage{amstext}
\usepackage{amssymb}
\usepackage{esint}
\usepackage{hyperref}

\makeatletter
\numberwithin{equation}{section}
\numberwithin{figure}{section}
\theoremstyle{plain}
\newtheorem{thm}{Theorem}[section]

\allowdisplaybreaks
\usepackage[sort&compress,square,numbers]{natbib}
\newcommand\relphantom[1]{\mathrel{\phantom{#1}}}

\makeatother

\global\long\def\Zp{\mathbb{Z}_{p}}

\begin{document}

\title[Symmetry for $q$-Bernoulli polynomials]{Some identities of symmetry for $q$-Bernoulli polynomials under
symmetric group of degree $n$}
\author{Dae San Kim}
\address{Department of Mathematics, Sogang University, Seoul 121-742, Republic
of Korea}
\email{dskim@sogang.ac.kr}

\author{Taekyun Kim}
\address{Department of Mathematics, Kwangwoon University, Seoul 139-701, Republic
of Korea}
\email{tkkim@kw.ac.kr}


\begin{abstract}
In this paper, we give some new identities of symmetry for $q$-Bernoulli
polynomials under the symmetric group of degree $n$ arising from
$p$-adic $q$-integrals on $\mathbb{Z}_{p}$.
\end{abstract}

\maketitle

\section{Introduction}

Let $p$ be a fixed prime number. Throughout this paper $\Zp$, $\mathbb{Q}_{p}$
and $\mathbb{C}_{p}$ will denote the ring of $p$-adic integers,
the field of $p$-adic rational numbers and the completion of the
algebraic closure of $\mathbb{Q}_{p}$. The $p$-adic norm is normalized
as $\left|p\right|_{p}=\frac{1}{p}$. Let $q$ be an indeterminate
in $\mathbb{C}_{p}$ such that $\left|1-q\right|_{p}<p^{-\frac{1}{p-1}}$.
The $q$-analogue of the number $x$ is defined as $\left[x\right]_{q}=\frac{1-q^{x}}{1-q}$.
Note that $\lim_{q\rightarrow1}\left[x\right]_{q}=x$. Let $UD\left(\Zp\right)$
be the space of uniformly differentiable functions on $\Zp$. For
$f\in UD\left(\Zp\right),$ the $p$-adic $q$-integral on $\Zp$
is defined by Kim as
\begin{equation}
I_{q}\left(f\right)=\int_{\Zp}f\left(x\right)d\mu_{q}\left(x\right)=\lim_{N\rightarrow\infty}\frac{1}{\left[p^{N}\right]_{q}}\sum_{x=0}^{p^{N}-1}f\left(x\right)q^{x},\quad\left(\text{see \cite{key-7}}\right).\label{eq:1}
\end{equation}

From (\ref{eq:1}), we have
\begin{equation}
qI_{q}\left(f_{1}\right)-I_{q}\left(f\right)=\left(q-1\right)f\left(0\right)+\frac{q-1}{\log q}f^{\prime}\left(0\right),\quad\text{where }f_{1}\left(x\right)=f\left(x+1\right).\label{eq:2}
\end{equation}

As is well known, the Bernoulli numbers are defined by
\[
B_{0}=1,\quad\left(B+1\right)^{n}-B_{n}=\begin{cases}
1 & \text{if }n=1,\\
0 & \text{if }n>1,
\end{cases}
\]
 with the usual convention about replacing $B^{n}$ by $B_{n}$ (see
\cite{key-1,key-2,key-3,key-4,key-5,key-6,key-7,key-8,key-9,key-10,key-11,key-12}).
The Bernoulli polynomials are given by
\begin{equation}
B_{n}\left(x\right)=\sum_{l=0}^{n}\binom{n}{l}B_{l}x^{n-l},\quad\left(n\ge0\right),\quad\left(\text{see \cite{key-11}}\right).\label{eq:3}
\end{equation}

In \cite{key-4}, L. Carlitz considered the $q$-analogue of Bernoulli
numbers as follows:
\begin{equation}
\beta_{0,q}=1,\quad q\left(q\beta_{q}+1\right)^{n}-\beta_{n,q}=\begin{cases}
1 & \text{if }n=1,\\
0 & \text{if }n>1,
\end{cases}\label{eq:4}
\end{equation}
with the usual convention about replacing $\beta_{q}^{n}$ by $\beta_{n,q}$.

He also defined $q$-Bernoulli polynomials as follows:
\begin{equation}
\beta_{n,q}\left(x\right)=\sum_{l=0}^{n}\binom{n}{l}q^{lx}\left[x\right]_{q}^{n-l}\beta_{l,q}\quad\left(\text{see \cite{key-2,key-3,key-4,key-8}}\right).\label{eq:5}
\end{equation}

In \cite{key-7}, Kim proved the following integral representation
related to Carlitz $q$-Bernoulli polynomials:
\begin{equation}
\beta_{n,q}\left(x\right)=\int_{\Zp}\left[x+y\right]_{q}^{n}d\mu_{q}\left(x\right),\quad\left(n\ge0\right).\label{eq:6}
\end{equation}

From (\ref{eq:2}), we note that
\begin{equation}
q\int_{\Zp}\left[x+1\right]_{q}^{n}d\mu_{q}\left(x\right)-\int_{\Zp}\left[x\right]_{q}^{n}d\mu_{q}\left(x\right)=\begin{cases}
q-1 & \text{if }n=0\\
1 & \text{if }n=1\\
0 & \text{if }n>1
\end{cases}.\label{eq:7}
\end{equation}

By (\ref{eq:7}), we get
\[
\beta_{0,q}=1,\quad q\beta_{n,q}\left(1\right)-\beta_{n,q}=\begin{cases}
1 & \text{if }n=1\\
0 & \text{if }n>1
\end{cases}.
\]

The purpose of this paper is to give identities of symmetry for Carlitz's
$q$-Bernoulli polynomials under the symmetric group of degree $n$
arising from $p$-adic $q$-integrals on $\Zp$.

\section{Symmetric identiites of $\beta_{n,q}\left(x\right)$ under $S_{n}$}

For $n\in\mathbb{N}$, let $w_{1},w_{2},\dots,w_{n}\in\mathbb{N}$.
Then, we have
\begin{align}
 & \relphantom{=}\int_{\Zp}e^{\left[\left(\prod_{j=1}^{n-1}w_{j}\right)y+\left(\prod_{j=1}^{n}w_{j}\right)x+w_{n}\sum_{j=1}^{n-1}\left(\prod_{\substack{i=1\\
i\neq j
}
}^{n-1}w_{i}\right)k_{j}\right]_{q}t}d\mu_{q^{w_{1}w_{2}\cdots w_{n-1}}}\left(y\right)\label{eq:8}\\
 &= \lim_{N\rightarrow\infty}\frac{1}{\left[p^{N}\right]_{q^{w_{1}w_{2}\cdots w_{n-1}}}}\sum_{y=0}^{p^{N}-1}e^{\left[\left(\prod_{j=1}^{n-1}w_{j}\right)y+\left(\prod_{j=1}^{n}w_{j}\right)x+w_{n}\sum_{j=1}^{n-1}\left(\prod_{\substack{i=1\\
i\neq j
}
}^{n-1}w_{i}\right)k_{j}\right]_{q}t}\nonumber \\
 &\relphantom{=} \times q^{\left(\prod_{j=1}^{n-1}w_{j}\right)y}.\nonumber \\
  &= \lim_{N\rightarrow\infty}\frac{1}{\left[w_{n}p^{N}\right]_{q^{w_{1}w_{2}\cdots w_{n-1}}}}\sum_{m=0}^{w_{n}-1}\sum_{y=0}^{p^{N}-1}\nonumber \\
 &\relphantom{=} \times e^{\left[\left(\prod_{j=1}^{n-1}w_{j}\right)\left(m+w_{n}y\right)+\left(\prod_{j=1}^{n}w_{j}\right)x+w_{n}\sum_{j=1}^{n-1}\left(\prod_{\substack{i=1\\
i\neq j
}
}^{n-1}w_{i}\right)k_{j}\right]_{q}t}q^{w_{1}w_{2}\cdots w_{n-1}\left(m+w_{n}y\right)}.\nonumber
\end{align}

Thus, by \ref{eq:8}, we get
\begin{align}
 &\relphantom{=} \frac{1}{\left[\prod_{l=1}^{n-1}w_{l}\right]_{q}}\prod_{l=1}^{n-1}\sum_{k_{l}=0}^{w_{l}-1}q^{w_{n}\sum_{j=1}^{n-1}\left(\prod_{\substack{i=1\\
i\neq j
}
}^{n-1}w_{i}\right)k_{j}}\label{eq:9}\\
&\times\int_{\Zp}e^{\left[\left(\prod_{j=1}^{n-1}w_{j}\right)y+\prod_{j=1}^{n}w_{j}x+w_{n}\sum_{j=1}^{n-1}\left(\prod_{\substack{i=1\\
i\neq j
}
}^{n-1}w_{i}\right)k_{j}\right]_{q}t}d\mu_{q^{w_{1}w_{2}\cdots w_{n-1}}}\left(y\right)\nonumber\\
 &= \lim_{N\rightarrow\infty}\frac{1}{\left[\prod_{l=1}^{n}w_{l}p^{N}\right]_{q}}\prod_{l=1}^{n-1}\sum_{k_{l}=0}^{w_{l}-1}\sum_{m=0}^{w_{n}-1}\sum_{y=0}^{p^{N}-1}q^{\left(\prod_{j=1}^{n-1}w_{j}\right)\left(m+w_{n}y\right)+\sum_{j=1}^{n-1}\left(\prod_{\substack{i=1\\
i\neq j
}
}^{n-1}w_{i}\right)k_{j}w_{n}}\nonumber \\
 & \times e^{\left[\left(\prod_{j=1}^{n-1}w_{j}\right)\left(m+w_{n}y\right)+\left(\prod_{j=1}^{n}w_{j}\right)x+w_{n}\sum_{j=1}^{n-1}\left(\prod_{\substack{i=1\\
i\neq j
}
}^{n-1}w_{i}\right)k_{j}\right]_{q}t}.\nonumber
\end{align}

We note that (\ref{eq:9}) is invariant under any permutation $\sigma\in S_{n}$.
Therefore, by (\ref{eq:9}), we obtain the following theorem.
\begin{thm}
\label{thm:1} For $w_{1},w_{2},\dots,w_{n}\in\mathbb{N}$, the following
expressions
\begin{align*}
 & \frac{1}{\left[\prod_{l=1}^{n-1}w_{\sigma\left(l\right)}\right]_{q}}\prod_{l=1}^{n-1}\sum_{k_{l}=0}^{w_{\sigma\left(l\right)}-1}q^{w_{\sigma\left(n\right)}\sum_{j=1}^{n-1}\left(\prod_{\substack{i=1\\
i\neq j
}
}^{n-1}w_{\sigma\left(i\right)}\right)k_{j}}\\
 & \times\int_{\Zp}e^{\left[\left(\prod_{j=1}^{n-1}w_{\sigma\left(j\right)}\right)y+\prod_{j=1}^{n}w_{j}x+w_{\sigma\left(n\right)}\sum_{j=1}^{n-1}\left(\prod_{\substack{i=1\\
i\neq j
}
}^{n-1}w_{\sigma\left(i\right)}\right)k_{j}\right]_{q}t}d\mu_{q^{w_{\sigma\left(1\right)}\cdots w_{\sigma\left(n-1\right)}}}\left(y\right)
\end{align*}
are the same for any $\sigma\in S_{n}$.
\end{thm}
We observe that
\begin{align}
 & \left[\left(\prod_{j=1}^{n-1}w_{j}\right)y+\left(\prod_{j=1}^{n}w_{j}\right)x+w_{n}\sum_{j=1}^{n-1}\left(\prod_{\substack{i=1\\
i\neq j
}
}^{n-1}w_{i}\right)k_{j}\right]_{q}\label{eq:10}\\
= & \left[\prod_{j=1}^{n-1}w_{j}\right]_{q}\left[y+w_{n}x+\frac{w_{n}}{w_{1}}k_{1}+\cdots+\frac{w_{n}}{w_{n-1}}k_{n-1}\right]_{q^{w_{1}\cdots w_{n-1}}}\nonumber \\
= & \left[\sum_{j=1}^{n-1}w_{j}\right]_{q}\left[y+w_{n}x+\sum_{j=1}^{n-1}\frac{w_{n}}{w_{j}}k_{j}\right]_{q^{w_{1}\cdots w_{n-1}}}.\nonumber
\end{align}

Thus, by (\ref{eq:10}), we get

\begin{align}
 & \int_{\Zp}e^{\left[\left(\prod_{j=1}^{n-1}w_{j}\right)y+\left(\prod_{j=1}^{n}w_{j}\right)x+w_{n}\sum_{j=1}^{n-1}\left(\prod_{\substack{i=1\\
i\neq j
}
}^{n-1}w_{i}\right)k_{j}\right]_{q}t}d\mu_{q^{w_{1}\cdots w_{n-1}}}\left(y\right)\label{eq:11}\\
= & \sum_{m=0}^{\infty}\left[\prod_{j=1}^{n-1}w_{j}\right]_{q}^{m}\int_{\Zp}\left[y+w_{n}x+w_{n}\sum_{j=1}^{n-1}\frac{k_{j}}{w_{j}}\right]_{q^{w_{1}\cdots w_{n-1}}}^{m}d\mu_{q^{w_{1}\cdots w_{n-1}}}\left(y\right)\frac{t^{m}}{m!}\nonumber \\
= & \sum_{m=0}^{\infty}\left[\prod_{j=1}^{n-1}w_{j}\right]_{q}^{m}\beta_{m,q^{w_{1}\cdots w_{n-1}}}\left(w_{n}x+\sum_{j=1}^{n-1}\frac{w_{n}}{w_{j}}k_{j}\right)\frac{t^{m}}{m!}.\nonumber
\end{align}

For $m\ge0$, from (\ref{eq:11}), we have
\begin{align}
 & \int_{\Zp}\left[\left(\prod_{j=1}^{n-1}w_{j}\right)y+\left(\prod_{j=1}^{n}w_{j}\right)x+w_{n}\sum_{j=1}^{n-1}\left(\prod_{\substack{i=1\\
i\neq j
}
}^{n-1}w_{i}\right)k_{j}\right]_{q}^{m}d\mu_{q^{w_{1}\cdots w_{n-1}}}\left(y\right)\label{eq:12}\\
= & \left[\prod_{j=1}^{n-1}w_{j}\right]_{q}^{m}\beta_{m,q^{w_{1}\cdots w_{n-1}}}\left(w_{n}x+\sum_{j=1}^{n-1}\frac{w_{n}}{w_{j}}k_{j}\right).\nonumber
\end{align}

Therefore, by Theorem \ref{thm:1} and (\ref{eq:12}), we obtain the
following theorem.
\begin{thm}
\label{thm:2} For $m\ge0$, $w_{1},\dots,w_{n}\in\mathbb{N}$, the
following expressions
\begin{align*}
&\left[\prod_{j=1}^{n-1}w_{\sigma\left(j\right)}\right]_{q}^{m-1}\prod_{l=1}^{n-1}\sum_{k_{l}=0}^{w_{\sigma\left(l\right)}-1}q^{\sum_{j=1}^{n-1}\left(\prod_{\substack{i=1\\
i\neq j
}
}^{n-1}w_{\sigma\left(i\right)}\right)k_{j}w_{\sigma\left(n\right)}}\\
&\times\beta_{m,q^{w_{\sigma\left(1\right)}\cdots w_{\sigma\left(n-1\right)}}}\left(w_{\sigma\left(n\right)}x+w_{\sigma\left(n\right)}\sum_{j=1}^{n-1}\frac{k_{j}}{w_{\sigma\left(j\right)}}\right)
\end{align*}
are the same for any $\sigma\in S_{n}$.
\end{thm}
It is easy to show that
\begin{align}
 & \left[y+w_{n}x+w_{n}\sum_{j=1}^{n-1}\frac{k_{j}}{w_{j}}\right]_{q^{w_{1}\cdots w_{n-1}}}\label{eq:13}\\
= & \frac{\left[w_{n}\right]_{q}}{\left[\prod_{j=1}^{n-1}w_{j}\right]_{q}}\left[\sum_{j=1}^{n-1}\left(\prod_{\substack{i=1\\
i\neq j
}
}^{n-1}w_{i}\right)k_{j}\right]_{q^{w_{n}}}\nonumber\\
&+q^{w_{n}\sum_{j=1}^{n-1}\left(\prod_{\substack{i=1\\
i\neq j
}
}^{n-1}w_{i}\right)k_{j}}\left[y+w_{n}x\right]_{q^{w_{1}\cdots w_{n-1}}}.\nonumber
\end{align}

From (\ref{eq:13}), we can derive the following equation:
\begin{align}
 &\relphantom{=} \int_{\Zp}\left[y+w_{n}x+w_{n}\sum_{j=1}^{n-1}\frac{k_{j}}{w_{j}}\right]_{q^{w_{1}\cdots w_{n-1}}}^{m}d\mu_{q^{w_{1}\cdots w_{n-1}}}\left(y\right)\label{eq:14}\\
 &= \sum_{l=0}^{m}\binom{m}{l}\left(\frac{\left[w_{n}\right]_{q}}{\left[\prod_{j=1}^{n-1}w_{j}\right]_{q}}\right)^{m-l}\left[\sum_{j=1}^{n-1}\left(\prod_{\substack{i=1\\
i\neq j
}
}^{n-1}w_{i}\right)k_{j}\right]_{q^{w_{n}}}^{m-l}q^{lw_{n}\sum_{j=1}^{n-1}\left(\prod_{\substack{i=1\\
i\neq j
}
}^{n-1}w_{i}\right)k_{j}}\nonumber \\
 &\relphantom{=} \times\int_{\Zp}\left[y+w_{n}x\right]_{q^{w_{1}\cdots w_{n-1}}}^{l}d\mu_{q^{w_{1}\cdots w_{n-1}}}\left(y\right)\nonumber \\
 &= \sum_{l=0}^{m}\binom{m}{l}\left(\frac{\left[w_{n}\right]_{q}}{\left[\prod_{j=1}^{n-1}w_{j}\right]_{q}}\right)^{m-l}\left[\sum_{j=1}^{n-1}\left(\prod_{\substack{i=1\\
i\neq j
}
}^{n-1}w_{i}\right)k_{j}\right]_{q^{w_{n}}}^{m-l}\\
&\relphantom{=}\times q^{lw_{n}\sum_{j=1}^{n-1}\left(\prod_{\substack{i=1\\
i\neq j
}
}^{n-1}w_{i}\right)k_{j}}\beta_{l,q^{w_{1}\cdots w_{n-1}}}\left(w_{n}x\right).\nonumber
\end{align}

Thus, by (\ref{eq:14}), we get
\begin{align*}
 &\relphantom{=} \left[\prod_{j=1}^{n-1}w_{j}\right]_{q}^{m-1}\prod_{l=1}^{n-1}\sum_{k_{l}=0}^{w_{l}-1}q^{w_{n}\sum_{j=1}^{n-1}\left(\prod_{\substack{i=1\\
i\neq j
}
}^{n-1}w_{i}\right)k_{j}}\\
&\relphantom{=}\times\int_{\Zp}\left[y+w_{n}x+w_{n}\sum_{j=1}^{n-1}\frac{k_{j}}{w_{j}}\right]_{q^{w_{1}\cdots w_{n-1}}}^{m}d\mu_{q^{w_{1}\cdots w_{n-1}}}\left(y\right)\\
&= \sum_{l=0}^{m}\binom{m}{l}\left[\prod_{j=1}^{n-1}w_{j}\right]_{q}^{l-1}\left[w_{n}\right]_{q}^{m-l}\beta_{l,q^{w_{1}\cdots w_{n-1}}}\left(w_{n}x\right)\nonumber\\
&\relphantom{=}\times\prod_{s=1}^{n-1}\sum_{k_{s}=0}^{w_{s}-1}q^{\sum_{j=1}^{n-1}\left(\prod_{\substack{i=1\\
i\neq j
}
}^{n-1}w_{i}\right)k_{j}w_{n}\left(l+1\right)}\left[\prod_{j=1}^{n-1}\left(\prod_{\substack{i=1\\
i\neq j
}
}^{n-1}w_{i}\right)k_{j}\right]_{q^{w_{n}}}^{m-l}\\
 &= \sum_{l=0}^{m}\binom{m}{l}\left[\prod_{j=1}^{n-1}w_{j}\right]_{q}^{l-1}\left[w_{n}\right]_{q}^{m-l}\beta_{l,q^{w_{1}\cdots w_{n-1}}}\left(w_{n}x\right)T_{m,q^{w_{n}}}\left(w_{1},w_{2},\dots,w_{n-1}\mid l\right),
\end{align*}
 where
\begin{align*}
 & T_{m,q}\left(w_{1},w_{2},\dots,w_{n-1}\mid l\right)\\
= & \prod_{s=1}^{n-1}\sum_{k_{s}=0}^{w_{s}-1}q^{\left(l+1\right)\sum_{j=1}^{n-1}\left(\prod_{\substack{i=1\\
i\neq j
}
}^{n-1}w_{i}\right)k_{j}}\left[\sum_{j=1}^{n-1}\left(\prod_{\substack{i=1\\
i\neq j
}
}^{n-1}w_{i}\right)k_{j}\right]_{q}^{m-l}.
\end{align*}
As this expression is invariant under any permutation in $S_{n}$, we have the following
theorem.
\begin{thm}
\label{thm:3} For $m\ge0$, $n,w_{1},\dots,w_{n}\in\mathbb{N}$,
the following expressions
\begin{align*}
&\sum_{l=0}^{m}\binom{m}{l}\left[\prod_{j=1}^{n-1}w_{\sigma\left(j\right)}\right]_{q}^{l-1}\left[w_{\sigma\left(n\right)}\right]_{q}^{m-l}\beta_{l,q^{w_{\sigma\left(1\right)}\cdots w_{\sigma\left(n-1\right)}}}\left(w_{\sigma\left(n\right)}x\right)\\
&\times T_{m,q^{w_{\sigma\left(n\right)}}}\left(w_{\sigma\left(1\right)},\dots,w_{\sigma\left(n-1\right)}\mid l\right)
\end{align*}
are all the same for $\sigma\in S_{n}$.

\end{thm}

\bibliographystyle{amsplain}

\providecommand{\bysame}{\leavevmode\hbox to3em{\hrulefill}\thinspace}
\providecommand{\MR}{\relax\ifhmode\unskip\space\fi MR }
% \MRhref is called by the amsart/book/proc definition of \MR.
\providecommand{\MRhref}[2]{%
  \href{http://www.ams.org/mathscinet-getitem?mr=#1}{#2}
}
\providecommand{\href}[2]{#2}

\end{document}